\documentstyle[amscd,12pt,psamsfonts]{amsart}
           \input xypic

\newtheorem{defn}{Definition}
\newtheorem{thm}[defn]{Theorem}

\newtheorem{cor}[defn]{Corollary}
\newtheorem{lem}[defn]{Lemma}
\newtheorem{prop}[defn]{Proposition}

\theoremstyle{remark}
\newtheorem{rem}{Remark}
\theoremstyle{remark}

\numberwithin{equation}{section}
\numberwithin{defn}{section}

\newcommand\Det{\operatorname{Det}}
\renewcommand\det{\operatorname{det}}

\renewcommand\deg{\operatorname{deg}}

\newcommand\supp{\operatorname{supp}}

\renewcommand\o{{\mathcal O}}
\renewcommand\P{{\mathcal P}}

\newcommand\M{{\mathcal M}}
\newcommand\U{{\mathcal U}}
\newcommand\SU{{\mathcal SU}}

\newcommand\C{{\mathbb C}}

\newcommand\wtilde{\widetilde}

\newcommand\beq{
                \setcounter{equation}{\value{defn}}\addtocounter{defn}1
                \begin{equation}}

\begin{document}

\title[Formulae for Non-Abelian Thetas]{Addition Formulae
for  Non-Abelian \\ Theta Functions and Applications}
\author[E. G\'omez Gonz\'alez and F. J. Plaza Mart\'{\i}n]{E. G\'omez
   Gonz\'alez  
   \\ F. J. Plaza Mart\'{\i}n  
}
   \address{Departamento de Matem\'aticas, Universidad de Salamanca,  Plaza
   de la Merced 1-4
   \\
   37008 Salamanca. Spain.
   \\
    Tel: +34 923294460. Fax: +34 923294583}

\thanks{
\indent This work is partially supported by the research contracts
BFM2000-1327 and  BFM2000-1315 of DGI and  SA064/01 of JCyL. The
second author is also supported by MCYT ``Ram\'on y Cajal''
program. \\
\indent This version will appear in Journal of Geometry and
Physics.}
\email{esteban@@usal.es}
\email{fplaza@@usal.es}

\begin{abstract}
This paper generalizes for non-abelian theta functions a number of
formulae valid for theta functions of Jacobian varieties. The
addition formula, the relation with the Sz\"ego kernel and with
the multicomponent KP hierarchy and the behavior under cyclic
coverings are given.
\end{abstract}

\maketitle



\section{Introduction}

Fay's addition formula for theta functions of Jacobians
(\cite{Fay} has turned out to be highly relevant in a number of
problems: geometric properties of Jacobians (existence of
trisecants to their Kummer varieties), infinitesimal behavior of
theta functions of Jacobians (KP and KdV equations) and algebraic
formulations of certain aspects of conformal field theories. On
the other hand, deep relations between moduli spaces of vector
bundles and Jacobians varieties (\cite{BNR,Li}) has been already
stablished. Therefore, it is thus natural to expect similarities
between the properties of classical theta functions and those of
non-abelian theta functions, in particular, an analogue of Fay's
addition formula.

In fact,  the existence of generalized addition formulae for
non-abelian theta functions has been conjectured by Schork
(conjecture IV.9 of \cite{Schork}) when generalizing for higher
rank Raina's approach to $b-c$ systems (\cite{Raina-bc}).
Therefore,  this kind of formulae should be useful tools when
studying Schork's correlation functions as well as non-abelian
generalizations of multiplicative Ward identities given by Witten
(\cite{Witten}). In fact, the case of line bundles has been
already worked out completely by Raina
(\cite{Raina-Corr,Raina-bc}). Further, in this direction, the
relations between theta functions and the Szeg\"o kernel are well
known (e.g. Theorem~25 of \cite{Hejhal} or \S6 of
\cite{Raina-Corr}) and have been useful in many problems (e.g.
\cite{Raina-bc}).

Moreover, an infinitesimal version of such a formula has been
given in (\cite{Pl}) when proving that non-abelian theta functions
verify the multicomponent KP hierarchy. Hence, an addition formula
may help in the understanding of this result and of its geometric
consequences (see \cite{MP} for the rank $1$ case).

On the other hand, the study of how Jacobian theta functions  vary
under morphisms of curves has shed light on their properties (e.g.
chapters IV and V of \cite{Fay}). This question is related to the
so-called twist structures of $b-c$ systems (\cite{Raina-bc}) and
is also addressed in pg. 844 of \cite{Ball} for higher rank.

The above problems are treated in this paper as follows. A
generalization of the addition formula for non-abelian theta
functions is the main result of section 3 (Theorem
\ref{thm:add-theta}) which coincides with Corollary~2.19 of
\cite{Fay} in the case of line bundles. This formula will be
derived as an identity among global sections of certain isomorphic
line bundles. It is worth mentioning some results needed for its
proof,  Theorem~\ref{thm:juugneneral} and
Theorem~\ref{thm:main-iso}, which allow us to determine the
pullback of the generalized theta divisor by different morphisms
which are essentially given by the action of the Jacobian on the
moduli space of semistable vector bundles. The latter theorem has
been already applied by Schork in~\cite{Schork2} in the study of
correlation functions of generalized $b-c$ systems.

The known relations between theta functions and the Szeg\"o kernel
associated to a line bundle are generalized in the fourth section.
The identity given in Theorem~\ref{thm:add-szego} could be a first
step in the question addressed by Ball and Vinnikov (pg. 865 of
\cite{Ball}) about the existence of a explicit formula for the
Szeg\"o kernel in higher rank. Another methods were used by Fay
(\cite{Fay2}) to give a similar relation for degree $0$ stable
vector bundles.

Theorem~\ref{thm:multicomponente} of the following section
contains a global version of Lemma~2.7 of \cite{Pl}. Recall that
the bilinear identity for the multicomponent KP hierarchy is a
consequence of this kind of formulae and that, in particular, the
non-abelian theta function is a tau-function of that hierarchy.

   The sixth section studies the behavior non-abelian theta functions
   under direct and inverse image by a cyclic covering
(Theorem~\ref{thm:cyclic-direct} and
Proposition~\ref{prop:inverseimage}).
   Since our methods are valid
for all $r\geq 1$, some of our results specialize to formulae of
the Jacobian case ($r=1$) given by Fay (see Remark
\ref{rem:prop5.1defay}).

{\small {\it Acknowledgments}: Both authors wish to express their
gratitude to Prof. J. M. Mu\~noz Porras for his valuable suggestions and
comments.}

\section{Preliminaries}

This section fixes notations and summarizes some results
   concerned with the generalized theta
divisor and non-abelian theta functions (see
\cite{BNR,DN,Laszlo,Po}).

Let $C$ be a irreducible smooth projective curve of genus $g\geq
2$ over $\C$. Given two integers, $r,d$, let $\U_C(r,d)$ (or
simply $\U(r,d)$) denote the moduli space of semistable vector
bundles on $C$ of rank $r$ and degree $d$. Let $\delta$ be
$p.g.c.d.(r,d)$ and $\bar r$ be $\frac{r^2}{\delta}$.

Recall that there is a closed subscheme of $\U(r,r(g-1))$ of codimension
1 given  by:
$$\Theta_r\,:=\, \{M\in \U(r,r(g-1))\,\colon\, h^0(C,M)> 0\}$$
It thus defines a
polarization which is called the generalized theta divisor (\cite{DN}).
Moreover, it holds that (Theorem~2 of \cite{BNR}):
$$h^0(\U(r,r(g-1)),\o(\Theta_r))\,=\, 1$$
We have therefore a global section $\theta_r$, defined up to a
constant, of $\o(\Theta_r)$ whose zero divisor is $\Theta_r$. This
section is known as the non-abelian theta function of rank $r$ and
degree $r(g-1)$.

   From \cite{DN,Po} we learn that in order to define a polarization on the
moduli space $\U(r,d)$ for an arbitrary $d$, we need to fix a
vector bundle $\bar F$ of degree
$\frac{-d}\delta+\frac{r}\delta(g-1)$  and rank $\frac{r}\delta$
such that there exists a vector bundle $E\in\U(r,d)$ with
$h^0(C,E\otimes\bar F)=0$. In particular, one obtains that
$\chi(M\otimes \bar F)=0$ for all $M\in\U(r,d)$.

  From now on, we will fix a theta characteristic $\o(\eta)$ on $C$
and we write $\bar F$ as $F(\eta)$ for a degree $-\frac{d}\delta$
rank $\frac{r}\delta$  vector bundle $F$. Then, it is known that:
$$\Theta_{[F_\eta ]}\,:=\, \{M\in \U(r,d)\,\colon\, h^0(C,M\otimes
F(\eta) )> 0\}$$ defines a polarization, that depends only on the
class of $F(\eta)$ in the Grothendieck group of coherent algebraic
sheaves on $C$. This divisor is known as the generalized theta
divisor on $\U(r,d)$ defined by $F(\eta)$.

Assume that a polarization $\Theta_{[F_\eta ]}$ in $\U(r,d)$ is
given. Recall that there exists $E$ such that $h^i(C,E\otimes
F(\eta))=0$ ($i=0,1$). Then,  by Lemma~2.5 of \cite{Po}, it
follows that $F$ is semistable.

Being $F$ semistable, we can define the following morphism:
\beq
\begin{aligned}
\label{eq:mrdenmr2} \beta_{F_{\eta}}:\U(r,d) &\,\longrightarrow \,
\U(\bar r,\bar r(g-1)) \\ M& \,\longmapsto\, M\otimes F(\eta)
\end{aligned}
\end{equation}
since the tensor product of semistable vector bundles is
semistable (see Theorem~3.1.4 of~\cite{ML}). It holds that
$\beta_{F_{\eta}}^{-1}(\Theta_{\bar r})=\Theta_{[F_\eta ]}$. Then,
define the non-abelian theta function $\theta_{[F_\eta ]}$ as the
image of $\theta_{\bar r}$  by the induced morphism:
   {\small  \beq\label{eq:rest-r}
H^0(\U(\bar r,\bar r(g-1)),\o(\Theta_{\bar r}))\,
\longrightarrow\, H^0(\U(r,d),\o(\Theta_{[F_\eta ]}))
\end{equation}}

However, the construction of these divisors as determinantal
subvarieties (\cite{DN,Laszlo}) turns out to be an essential tool
when proving statements about them.

Let $S$ be a scheme and $M$ be a semistable vector bundle on
$C\times S$ of rank $r$ and degree $d$ and let $\phi_M$ be the
morphism:
$$
\begin{aligned}
\phi_M:S &\,\longrightarrow \, \U(r,d) \\
s& \,\longmapsto\, M\vert_{C\times\{s\}}
\end{aligned}
$$
Then, the polarization satisfies the following property:
$$\phi^*_M(\o(\Theta_{[F_\eta ]}))\,\simeq \,\Det\big(R^\bullet
q_*(M\otimes p^*( F(\eta) )\big)^*$$ where $q:C\times S\to S$ and
$p:C\times S\to C$ are the natural projections.

In order to compute this determinant we proceed as follows. Fix an
effective divisor $D$ on $C\times S$ such that $R^1q_*( M(D)
\otimes p^*F(\eta))=0$. Then, tensor with $M\otimes p^*F(\eta)$
the following exact sequence on $C\times S $:
$$0\to \o\to \o(D)\to \o_D(D)\to 0$$
and consider the induced cohomology sequence on $S$:
{\small
$$\begin{aligned} 0\to q_*(M\otimes p^*F(\eta)) \to q_*(M(D)\otimes
p^*F(\eta))\overset{\alpha}\to
   & \\
\to q_*(M\otimes\o_D(D)\otimes p^*F(\eta)) & \to R^1q_*(M\otimes
p^*F(\eta))\to 0 \end{aligned}$$}

The properties of determinants (\cite{KM}) show that (up to a
constant): {\small\beq \label{eq:thetadet} \phi_M^*(\theta_{F_\eta
})=\det(\alpha)\in H^0(S,\phi_M^*\o(\Theta_{[F_\eta ]}))
\end{equation}} which is an effective way to deal with non-abelian theta
functions.

Finally, it is worth pointing out that the above construction also
applies to the universal bundle of $\U(r,d)$ when  $r,d$ are
coprime.


\section{Addition Formula}

The first part of this section is devoted to the explicit
computation of the pullback of the generalized theta divisor
$\Theta_{[F_\eta]}$ by the natural morphism:
\beq\label{eq:tensorju}
\begin{aligned}
m: \U(r,d)\times J &\,\longrightarrow\, \U(r,d) \\
(M,L)&\,\longmapsto\, M\otimes L
\end{aligned}
\end{equation}
where $J$ denotes the Jacobian variety of $C$, that is, isomorphism
classes of degree $0$ line bundles.

This calculation requires a number of intermediate results. Let us
introduce the following notation. Let $J_d$ denote the variety of
isomorphism classes of degree $d$ line bundles on $C$. The choice
of the theta characteristic $\eta$ gives rise to a principal
polarization on $J$, $\Theta_J$. Denote by $\phi_{\Theta_J}:J\to
{\operatorname{Pic}}^0(J)$ the isomorphism induced by $\Theta_J$.

Consider the morphism:
$$\det\colon\U(r,r(g-1))\,\to\, J_{r(g-1)}$$ which maps a vector
bundle to its determinant. Finally, for a line bundle $L\in J$ let
$T_L$ denote the translation defined by $L$ on the moduli space of
vector bundles as well as on the Jacobian variety.

\begin{lem}\label{lem:polarization}
Let $L\in J$. Then, there is an isomorphism:
$$T_L^*(\o_{\U}(\Theta_r))\otimes\o_{\U}(-\Theta_r)\simeq \det^*(T_{-r\eta}^*(
\phi_{\Theta_J}(L)))$$
\end{lem}

\begin{pf}
Let $\SU (r,\o(r\eta))$ be the moduli space of semistable vector
bundles of rank $r$ with determinant isomorphic to $\o(r\eta)$ and
let $\bar\Theta_r$ be the restriction of $\Theta_r$ to
$\SU(r,\o(r\eta))$.

Since ${\operatorname{Pic}}(\U(r,r(g-1)))\simeq
{\operatorname{Pic}}(\SU(r,\o(r\eta)))\oplus
{\operatorname{Pic}}(J_{r(g-1)})$ and $
{\operatorname{Pic}}(\SU(r,\o(r\eta)))\simeq {\mathbb
Z}\bar\Theta_r$ (\cite{DN}), one has that
$T_L^*(\o_{\U}(\Theta_r))\otimes\o_{\U}(-\Theta_r)\simeq
\det^*(N)$ for some $N\in {\operatorname{Pic}}^0(J_{r(g-1)})$
depending on $L$.

   Consider the
morphism:
    $$
\begin{aligned}
\bar m: \SU(r,\o(r\eta))\times J &\,\longrightarrow\, \U(r,r(g-1)) \\
(M,L)&\,\longmapsto\, M\otimes L
\end{aligned}
$$
By \cite{BNR}, we know that $\bar
m^*(\o_{\U}(\Theta_{r}))\,\simeq\,
    p_{\SU}^*(\o_{\SU}(\bar\Theta_r))\otimes
p_J^*(\o_J(r\Theta_J))$, where $p_{\SU}$ and $p_J$ are the natural
projections.

Taking the pull-back of
$T_L^*(\o_{\U}(\Theta_r))\otimes\o_{\U}(-\Theta_r)$ by the map
$\bar m$, one checks that $T_{r\eta}^*(N)\simeq
\phi_{\Theta_J}(L)\otimes \mu$ with $\mu$ a $r$-torsion point of
${\operatorname{Pic}}^0(J)$ depending on $L$.  Since $J$ is
complete and the $r$-torsion subgroup of
${\operatorname{Pic}}^0(J)$ is finite, one obtains that $\mu$ does
not depend on $L$. Letting $L=\o_C$, one has that $\mu=\o_J$ and
the claim follows.
\end{pf}

Now, we consider the morphism \ref{eq:tensorju} for the case
$d=r(g-1)$. Recall that  the Poincar\'e bundle on $J\times J$ is
given by:
   $$\P\,:=\, m_J^*(\o_J(\Theta_J))\otimes
p_1^*(\o_J(-\Theta_J))\otimes p_2^*(\o_J(-\Theta_J))$$ where
$m_J:J\times J\to J$ corresponds to the tensor product of line
bundles and $p_i$ is the projection onto the $i$-th factor.

\begin{lem}\label{lem:juu}
It holds that: {\small $$m^*(\o_{\U}(\Theta_r))\,=\,
    p_{\U}^*(\o_{\U}(\Theta_r))\otimes p_J^*(\o_J(r\Theta_J)) \otimes
( (T_{-r\eta}\circ\det)\times 1)^*\P$$}
\end{lem}

\begin{pf}
We consider the bundle on $\U(r,r(g-1))\times J$ defined by:
$${\mathcal F}\,:=\, m^*\o_{\U}(\Theta_r)\otimes
p_{\U}^*(\o_{\U}(-\Theta_r))\otimes ( (T_{-r\eta}\circ\det)\times
1)^*\P^{-1}$$ Then,  the above lemma implies that: {\small
$$
{\mathcal F}\vert_{\U(r,r(g-1))\times\{L\}}\,=  \,
T_L^*(\o_{\U}(\Theta_r))\otimes \o_{\U}(-\Theta_r)\otimes
\det^*(T_{-r\eta}^*(\phi_{\Theta_J}(L)))^* \, \simeq \, \o_{\U}
$$}
where $L$ is a  point of $J$.

   Hence, ${\mathcal F}$ is
trivial along the fibres of the natural projection $p_J:
\U(r,r(g-1))\times J\to J$. Seesaw theorem allows  us to conclude
that ${\mathcal F}\simeq p_J^*N$ for some
$N\in{\operatorname{Pic}}(J)$.

If we show that $N\simeq\o_J(r\Theta_J)$, we are done. Recall from
\cite{Raynaud} that there exists a vector bundle $M\in
\U(r,r(g-1))$ with $\wedge M:=\det(M)=\o(r\eta)$ such that the
subscheme $D(M):= \{L\in J\,\colon\, h^0(M\otimes L)>0\}$ is a
divisor of $J$ which is linearly equivalent to $r\Theta_J$. We now
have that:
$$N\,\simeq\, {\mathcal F}\vert_{\{M\}\times J}\, =\, \o_J(D(M))\otimes
\P^{-1}\vert_{\{(\wedge M)\otimes\o(-r\eta)\}\times J}\, \simeq\,
\o_J(r\Theta_J)$$
\end{pf}

We are now ready to compute the pullback of the generalized theta
divisor by the morphism \ref{eq:tensorju}.

\begin{thm}\label{thm:juugneneral}
One has that: {\small $$m^*(\o_{\U}(\Theta_{[F_\eta]}))\,=\,
    p_{\U}^*(\o_{\U}(\Theta_{[F_\eta]}))\otimes p_J^*(\o_J(\bar  r\Theta_J))
\otimes ( (\det\circ\beta_F)\times 1)^*\P$$} where $\beta_F\colon
\U(r,d)\to \U(\bar r,0)$ corresponds to tensor product by $F$.
\end{thm}

\begin{pf}
It follows from the above lemma and the following commutative
diagram:
$$\xymatrix{ \U(r,d)\times J \ar[r]^-{m} \ar[d]_-{\beta_{F_\eta}\times 1}
& \U(r,d) \ar[d]^-{\beta_{F_\eta}}
\\ \U(\bar r,\bar r(g-1))\times J \ar[r]^-{m}  & \U(\bar r,\bar r(g-1))}$$
\end{pf}

\begin{cor}\label{cor:beta}
Let $M\in \U(r,d)$ and $\beta_M :J\to \U(r,d)$ be the morphism
which sends $L$ to $M\otimes L$. It holds  that:
$$\beta_M^*(\o_{\U}(\Theta_{[F_\eta]}))\simeq \o_J(\bar r \Theta_J)\otimes
\phi_{\Theta_J}(\wedge(M\otimes F))$$
\end{cor}

\begin{pf}
It follows from the previous theorem and from the isomorphism
$\beta_M^*(\o_{\U}(\Theta_{[F_\eta]})) \simeq
m^*(\o_{\U}(\Theta_{[F_\eta]}))\vert_{\{M\}\times J}$.
\end{pf}

The rest of this section aims at giving explicit formulas for the
pullback of non-abelian theta functions by the morphism:
$$\begin{aligned}
\alpha_M\, : \, C^{2m}
\,\longrightarrow& \, \U(r,d) \\
(x_1,y_1,\dots,x_m,y_m)\mapsto & M(\sum_{i=1}^m (x_i -y_i))
\end{aligned}$$
where $C^{2m}$ is the product of $2m$ copies of the curve $C$ and
$M\in\U(r,d)$.

   Firstly, we will deal
with an isomorphism of line bundles on $C^{2m}$ and then it will
be applied it to obtain an identity among global sections of them.
Such a formula can be understood as an addition formula for
non-abelian theta functions. For the rank $1$ case and identifying
the theta function (as a section) with its classical analytic
expression, our formula turns out to coincide with Fay's formula.
However, as long as no analytic expressions for non-abelian theta
functions are known, our generalization must be regarded as an
identity of sections.

If a point of $C^{2m}$ is denoted by $(x_1,y_1,\dots,x_m,y_m)\in
C^{2m}$, we will call an index odd (resp. even) if  it corresponds
to a variable $x_i$ (resp. $y_j$). Finally, let $p_i$ be the
projection onto the $i$-th factor and $\Delta_{ij}$ the divisor of
$C^{2m}$ where the $i$-th and the $j$-th entries coincide.

\begin{lem}\label{lem:raina}
Let $L\in J$ and consider the following morphism:
$$\begin{aligned}
\alpha_L\,:\,  C^{2m}\, \longrightarrow & \, J \\
(x_1,y_1,\dots,x_m,y_m)\mapsto & L(\sum_i (x_i-y_i))
\end{aligned}
$$

Then, one has that:
$$\begin{aligned}
\alpha_L^*\o_J(\Theta_J)\,\simeq\, &
\o(\sum_{i+j=\text{odd}}\Delta_{ij}-
\sum_{i+j=\text{even}}\Delta_{ij}) \otimes  \\ & \otimes
\big(\underset{i\text{ odd}}\otimes p_i^* L^*(\eta) \big)
\otimes\big(\underset{i\text{ even}}\otimes  p_i^* L(\eta)\big)
\end{aligned}
$$
where the sums involve only $i,j$ with $i<j$.
\end{lem}

\begin{pf}
See  Theorem~11.1 in \cite{Raina-Corr}.
\end{pf}

\begin{thm}\label{thm:main-iso}
Let $M$ be a rational point of $\U(r,d)$ and let $\alpha_M$ be the
morphism defined by:
$$\begin{aligned}
\alpha_M\, : \, C^{2m}
\,\longrightarrow& \, \U(r,d) \\
(x_1,y_1,\dots ,x_m,y_m)\mapsto & M(\sum_{i=1}^m (x_i -y_i))
\end{aligned}$$

Then, there is an isomorphism of line bundles on $ C^{2m}$:
{\small
$$\begin{aligned}
\alpha_M^*\o_{\U}(\Theta_{[F_\eta]})\,\simeq\,
     &
\o(\sum_{i+j=\text{odd}}\Delta_{ij}-
\sum_{i+j=\text{even}}\Delta_{ij})^{\otimes \bar r} \otimes  \\ &
\otimes \big(\underset{i\text{ odd}}\otimes p_i^*\wedge (M\otimes
F (-\eta))^*\big) \otimes\big(\underset{i\text{ even}}\otimes
p_i^*\wedge (M\otimes F (\eta))\big)
\end{aligned}
$$}
where the sums involves only $i,j$ with $i<j$.
\end{thm}

\begin{pf}
The morphism $\alpha_M$ factors as follows: {\small
$$ C^{2m} \overset{\alpha_{\o}}\longrightarrow J\simeq \{M\} \times J
\overset{\beta_M}\longrightarrow \U(r,d)$$}

Let $M'$ be $M\otimes F$. Recalling Corollary~\ref{cor:beta} and
Lemma~\ref{lem:raina} we have that: {\small
$$
\begin{aligned}
\alpha_M^*\o_{\U}(\Theta_{[F_\eta]})\,\simeq &\,\alpha_{\o}^*\big(
\o_J(\bar r
\Theta_J)\otimes \phi_{\Theta_J}(\wedge M')\big)\, \simeq\\
\simeq &\,
    \alpha_{\o}^*\big( T_{\wedge
M'}^*\o_J(\Theta_J)\otimes \o_J(\Theta_J)^{\otimes\bar r-1}\big) \, \simeq \\
\simeq &\, \alpha_{\wedge M'}^*(\o_J(\Theta_J))\otimes
\alpha_{\o_C}^*(\o_J(\Theta_J))^{\otimes\bar r-1} \, \simeq \\
\simeq &\,
\o(\sum_{i+j=\text{odd}}\Delta_{ij}-
\sum_{i+j=\text{even}}\Delta_{ij}) \otimes  \\
& \otimes \big(\underset{i\text{ odd}}\otimes p_i^*(\wedge
M')^*(\eta)\big) \otimes\big(\underset{i\text{
even}}\otimes  p_i^*(\wedge M')(\eta)\big) \otimes  \\
   & \otimes
\o(\sum_{i+j=\text{odd}}\Delta_{ij}-
\sum_{i+j=\text{even}}\Delta_{ij})^{\otimes\bar r-1} \otimes  \\
& \otimes \big(\underset{i\text{ odd}}\otimes
p_i^*\o(\eta)^{\otimes\bar r-1}\big) \otimes\big(\underset{i\text{
even}}\otimes p_i^*\o(\eta)^{\otimes\bar r-1}\big)
\end{aligned}$$}
And the theorem follows.
\end{pf}

\begin{rem}
This result has been applied in~\cite{Schork2} when proving the
relation of determinants of correlation functions of generalized
$b-c$ systems and determinants of non-abelian theta functions.
This is a first step of the expected fact that correlation
functions of generalized $b-c$-system are determined completely by
the geometry of the non-abelian theta divisor, analogously to the
rank one case.
\end{rem}

Recall from chapter~II of \cite{Fay} that the line bundle
$\o(\Delta)$ con $C\times C$ has a unique section $E(x,y)$, which
is known as the prime form and that it can be constructed in terms
of $\eta$. To be consistent with Fay, it will be assumed that the
theta characteristic $\eta$ is odd. In particular, it holds that
$E(x,y)=-E(y,x)$.

\begin{thm}\label{thm:add-theta}
Let $M$ be a rational point of $\U(r,d)$ such that $\theta_{F_\eta
}(M)\neq 0$. Then, for $(x_1,y_1,\dots,x_m,y_m)\in C^{2m}$, one
has that: {\small $$\begin{aligned}
\frac{\theta_{F_\eta }(M(\sum_{i=1}^m (x_i-y_i)
))}{\theta_{F_\eta }(M)}\cdot &
\prod_{i<j}E(x_i,x_j)^{\bar r}E(y_i,y_j)^{\bar r}
\,=\,  \\ \,=\, &
\prod_{i,j}E(x_i,y_j)^{\bar r} \cdot
\operatorname{det} \left(
\frac{\theta_{F_\eta }(M(x_i-y_j
))}{\theta_{F_\eta }(M)E(x_i,y_j)^{\bar r} }\right)
\end{aligned}$$}
\end{thm}

\begin{pf}
Observe that the r.h.s. of the equality of the statement equals
the sum:
$$ \sum_{\sigma\in{\mathcal S}_m}
\operatorname{sign}{(\sigma)} \prod \Sb i,j \\ \sigma(i)\neq j
\endSb  E(x_i,y_j)^{\bar r} \cdot \prod_i
\frac{\theta_{F_{\eta}}(M(x_i-y_{\sigma(i)} ))}{\theta_{F_{\eta}}
(M)}$$ By Theorem \ref{thm:main-iso} with $m=1$,
$\frac{\theta_{F_\eta } (M(x_i-y_j ))}{\theta_{F_\eta } (M)}$ is a
section of the line bundle: {\small
$$\o(\Delta)^{\otimes \bar r} \otimes p_i^*\wedge(M\otimes F(-\eta)
)^* \otimes p_j^*\wedge(M\otimes
F(\eta))$$} So, it turns out that each term of the above sum is a
global section of: {\small
$$
\o(\sum_{i+j=\text{odd}}\Delta_{ij})^{\otimes\bar r} \otimes
\big(\underset{i\text{ odd}}\otimes p_i^*\wedge (M\otimes F(-\eta)
)^*\big) \otimes\big(\underset{i\text{ even}}\otimes p_i^*\wedge
(M\otimes F (\eta) )\big)
$$}

The l.h.s. is a section of the line bundle:
{\small
$$
\alpha_M^*\o(\Theta_{[F_\eta ]})\otimes
\o(\sum_{i+j=\text{even}}\Delta_{ij})^{\otimes\bar r}
$$}

These two line bundles are isomorphic by Theorem \ref{thm:main-iso}. Hence,
both sides of the equality are to be understood as global sections of the same
line bundle. Since $C^{2m}$ is proper, there is no non-constant global
section of the trivial bundle. So, if we show that both
sections have the same zero divisor, then they coincide up to a
multiplicative constant, which will be eventually shown to be $1$.

Let $D_R$ (resp. $D_L$) denote the zero divisor of the r.h.s.
(resp. l.h.s.)~. Since $D_L$ and $D_R$ are linearly equivalent,
there exists a rational function $f$ on $C^{2m}$ such that:
$$D_R-D_L\,=\,D(f)$$

Let us consider the following diagram:
$$\xymatrix{ C^{2m} \ar[r]^-{f} \ar[d]_-{\pi} & {\mathbb P}^1 \\
C^{2m-1}  &}$$ where
$\pi(x_1,y_1,\dots,x_m,y_m):=(x_1,x_2,y_2,\dots,x_m,y_m)$.

Suppose we have proved that there exists $z\in C^{2m-1}$ such that
$D_L\vert_{\pi^{-1}(z)}= D_R\vert_{\pi^{-1}(z)}$ and
$\operatorname {supp}( D_L\vert_{\pi^{-1}(z)})\neq \pi^{-1}(z)$.
It thus follows that $f\vert_{\pi^{-1}(z)}$ is a non-zero
constant, since ${\pi^{-1}(z)}\simeq C$ is proper. From the
Rigidity Lemma one deduces that $f$ is constant  along the fibers
of $\pi$ and, therefore, $f$ has not poles nor zeroes in $C^{2m}$.
Summing up, $f$ is invertible or, what amounts to the same, $D_L=
D_R$. So, there exists a non-zero constant $\lambda$ such that the
l.h.s. equals the r.h.s. multiplied  by $\lambda$. Letting
$x_i=y_i$ for all $i$, we obtain that $\lambda=1$.

By the above discussion, it remains to show that there exists $z$
such that $D_L\vert_{\pi^{-1}(z)}= D_R\vert_{\pi^{-1}(z)}$ and
$\operatorname {supp} (D_L\vert_{\pi^{-1}(z)})\neq \pi^{-1}(z)$.

       We take
$z=(x_1,x_2,y_2,\dots,x_m,y_m)\in C^{2m-1}$ such that $x_k=y_k$
for $k\neq 1$ and $x_i\neq x_j$ for all $i\neq j$. Then, we have
that: {\small $$\begin{aligned} &\Big (\frac{\theta_{F_{\eta}}
(M(\sum_{i=1}^m (x_i-y_i) ))}{\theta_{F_{\eta}} (M)}\cdot
\prod_{i<j}E(x_i,x_j)^{\bar r}E(y_i,y_j)^{\bar r}\Big )
\Big\vert_{\pi^{-1}(z)}\, = \\ &\qquad =\, \frac{\theta_{F_{\eta}}
(M( x_1-y_1))}{\theta_{F_{\eta}} (M)}\cdot \prod_{k\neq
1}E(y_1,y_k)^{\bar r} \prod_{i<j}E(x_i,x_j)^{\bar r}\prod_{2\leq
i<j}E(y_i,y_j)^{\bar r}\end {aligned}$$} and the r.h.s. restricted
to the fibre  of $z$ is:
{\small $$\begin{aligned} \det  \Big ( \prod_{k\neq i}
E(x_k,y_j)^{\bar r}\cdot &\frac {\theta_{F_{\eta}}
(M(x_i-y_j))}{\theta_{F_{\eta}} (M)}\Big
)\Big\vert_{\pi^{-1}(z)}\, = \\ &\qquad\qquad =\,  \prod_{k\neq
1}E(x_k,y_1)^{\bar r}\cdot \frac {\theta_{F_{\eta}}
(M(x_1-y_1))}{\theta_{F_{\eta}} (M)}\end{aligned}$$}

Letting $y_1=x_1$ one checks that both restrictions are not zero.
Furthermore, since the first one is equal to the second times a
non-zero constant on $\pi^{-1}(z)$, one has that
$D_L\vert_{\pi^{-1}(z)}= D_R\vert_{\pi^{-1}(z)}$. The theorem is
proved.
\end{pf}

\section{Addition formula and the Szeg\"o kernel}

Now,  let us recall briefly the definition and properties of the
Szeg\"o kernel  associated to a vector bundle $M\in
\U(r,d)-\Theta_{[F_\eta]}$. For such a bundle define the Szeg\"o
kernel, $S_M(x,y)$, to be the meromorphic section of
$p_1^*(M\otimes F(-\eta))^*\otimes p_2^*(M\otimes F(\eta))$ on
$C\times C$ with a simple pole along the diagonal such that its
residue along it is $1$.

Note that $S_M(x,y)$ might be written as an $r\times r$ matrix,
because there is an isomorphism: {\small $$\begin{aligned}
p_1^*(M\otimes F(-\eta))^* & \otimes p_2^*(M\otimes F(\eta))
\,\simeq\\
& \simeq\, {\mathcal H}\!\textit{om}(p_1^*(M\otimes F(-\eta)),
p_2^*(M\otimes F(\eta)))
\end{aligned}$$}

On the other hand, observe that the restriction to the diagonal
$\Delta\subset C\times C$ induces an isomorphism: {\small
$$H^0(C\times C,p_1^*(M\otimes F(-\eta))^*\otimes p_2^*(M\otimes
F(\eta))\otimes\o(\Delta))\,\simeq \, H^0(C, {\mathcal
E}\!\textit{nd}(M\otimes F))$$} and denote by $S_M^h(x,y)$ the
holomorphic global section of the vector bundle $p_1^*(M\otimes
F(-\eta))^*\otimes p_2^*(M\otimes F(\eta))\otimes\o(\Delta)$ whose
image by the above isomorphism is the identity.

Then, it is worth noting that $ E(x,y)\cdot S_M(x,y)$ is a holomorphic
     section of $p_1^*(M\otimes F(-\eta))^*\otimes
p_2^*(M\otimes F(\eta))\otimes\o(\Delta)$, because the morphism
$\o\to\o(\Delta)$ maps the global section $1$ to the global
section $E(x,y)$. One checks that $S_M^h(x,y)- E(x,y)\cdot
S_M(x,y)$ gives a global section of $p_1^*(M\otimes
F(-\eta))^*\otimes p_2^*(M\otimes F(\eta))$. Since this bundle has
no non-zero section,  one then has that: $$ S_M^h(x,y)\,=\,
E(x,y)\cdot S_M(x,y)
$$
If $S^h_M$ and $S_M$ are both understood as matrices, then this
identity makes sense too.

\begin{rem}
One can show that the rows of $E(x_0,y)\cdot S_M(x_0,y)$ for a
fixed point $x_0\in C$ give a basis of $H^0(C,M\otimes
F(\eta+x_0))$, because the restriction to $\{x_j\}\times C$ maps
$S^h_M$ to its rows:
{\small $$\begin{aligned}  H^0\big( C\times C,
{\mathcal H}\!\textit{om}\; (p_1^*(M\otimes F(-\eta))\, ,\,  p_2^*
(M\otimes & F(\eta)) \otimes\o(\Delta))\big) \,\longrightarrow \\
   &\longrightarrow \,
     H^0(C,M\otimes F(\eta+ x_j))^{\oplus r}
     \end{aligned}$$}
\end{rem}

Now, the relation of the non-abelian theta function and the
Szeg\"o kernel given by Fay for degree $0$ stable bundles
(\cite{Fay2}) is generalized for semistable ones by the following:

\begin{thm}\label{thm:add-szego}
Let $M$ be a rational point of $\U(r,d)$ such that $\theta_{F_\eta
}(M)\neq 0$. Then, for $(x_1,y_1,\dots,x_m,y_m)\in C^{2m}$, one
has that: {\small $$
\begin{aligned}
\frac{\theta_{F_\eta }(M(\sum_{i=1}^m (x_i-y_i)
))}{\theta_{F_\eta }(M)}\cdot
\prod_{i<j}E(x_i,x_j)^{\bar r} & E(y_i,y_j)^{\bar r} \,=
\\ \,=
\,
&
\prod_{i,j}E(x_i,y_j)^{\bar r} \cdot
\operatorname{det} S_M(x,y)
\end{aligned}
$$}
where $S_M(x,y)$ is a ${\bar r}m\times {\bar r}m$ matrix builded up from
the
${\bar r}\times {\bar r}$ boxes $S_M(x_i,y_j)$.
\end{thm}

Before giving the proof we need some results.

\begin{lem}\label{lem:sz2}
Both sides of the equality in the statement of
Theorem~\ref{thm:add-szego} are global sections of the isomorphic
line bundles on $C^{2m}$.
\end{lem}

\begin{pf}
Note that the matrix $S_M(x,y)$ is a meromorphic section of the bundle:
$$
{\mathcal H}\!\textit{om}\;  \big( \underset{i=\text{odd}}\oplus
p_i^*(M\otimes F(-\eta))\, ,\, \underset{j=\text{even}}\oplus
p_j^*(M\otimes F(\eta))\big)$$
       with poles along
$\sum_{i+j=\text{odd}}\Delta_{ij}$ (odd indexes correspond to
$x$'s variables, while even indexes correspond to $y$'s
variables). Therefore, the determinant $\det S_M(x,y)$ is a
meromorphic section of: $$\big(\otimes_{i=\text{odd}} p_i^*\wedge
(M\otimes
F(-\eta))^*\big)\otimes\big(\otimes_{j=\text{even}}p_j^*\wedge
(M\otimes F(\eta))\big)$$
   Counting the order of these poles, one
concludes that the r.h.s. is a holomorphic section of:
$$\big(\otimes_{i=\text{odd}} p_i^*\wedge
(M\otimes
F(-\eta))^*\big)\otimes\big(\otimes_{j=\text{even}}p_j^*\wedge
(M\otimes
F(\eta))\big)\otimes\o(\sum_{i+j=\text{odd}}\Delta_{ij})^{\otimes\bar
r}
$$

The l.h.s. is a holomorphic global section of:
$$\alpha^*(\Theta_{[F_\eta]})\otimes
\o(\sum_{i+j=\text{even}}\Delta_{ij})^{\otimes\bar r} $$
and the two  line bundles above are isomorphic by
Theorem~\ref{thm:main-iso}.
\end{pf}

\begin{lem}\label{lem:sz3}
Let $M\in \U(r,d)-\Theta_{[F_\eta]}$. Then, for $(x,y)\in C\times
C$, one has that: {\small $$ \frac{\theta_{F_\eta
}(M(x-y))}{\theta_{F_\eta }(M)} \,= \, E(x,y)^{\bar r} \cdot
\operatorname{det} S_M(x,y)
$$}
\end{lem}

\begin{pf}
Lemma~\ref{lem:sz2} implies that both sides are global sections of
the same line bundle.

Label the three copies of $C$ in $C\times C\times C$ by $0,1$ and
$2$. Let $\Delta_{0i}$ be the subscheme where the $0$-th entry
coincides with the $i$-th entry.  Let $p$ denote the projection
from  $C\times C\times C$ onto the copy of $C$ labelled with $0$.
Finally, let $q:C\times C\times C\to C\times C$ be the projection
onto the copies labelled with $1$ and $2$.

The bundle $M$ defines the morphism:
$$\begin{aligned}
\alpha_M:C\times C &\,\longrightarrow\,\U(r,d) \\
(x,y)&\,\longmapsto M(x-y)
\end{aligned}$$
By the construction of the polarization it is known that:
$$\alpha_M^*(\o(-\Theta_{[F_\eta]}))\,\simeq\, \Det(R^\bullet
q_*{\mathcal M})$$ where ${\mathcal M}:=p^*(M\otimes
F(\eta))\otimes\o(\Delta_{01}-\Delta_{02})$. Let $M'$ be $M\otimes
F(\eta)$.

Let us compute this determinant as well as a section of its dual. Consider
the exact sequence on $C\times C\times C$:
$$0\to \o(\Delta_{01}-\Delta_{02})\to \o(\Delta_{01})\to
\o(\Delta_{01})\vert_{\Delta_{02}}\to 0$$ Tensor with $p^*M'$ and
pushing it forward by $q$ one obtains: {\small $$ 0\to
q_*{\mathcal M}\to q_*(p^*M'\otimes \o(\Delta_{01}))
\overset{\beta}\to  q_*((p^*M'\otimes
\o(\Delta_{01})\vert_{\Delta_{02}})  \to  R^1q_*{\mathcal M}\to 0
$$}
   because $R^1 q_*(p^*M'\otimes \o(\Delta_{01})=0$.

Observe that the two middle terms of the above sequence are
locally free of same rank. Then, it follows that there exists a
canonical isomorphism: {\small
$$\alpha_M^*(\o(\Theta_{[F_\eta]}))\,\simeq\, \wedge\big(
q_*(p^*M'\otimes \o(\Delta_{01}))\big)^* \otimes \wedge
q_*(p^*M'\otimes \o(\Delta_{01})\vert_{\Delta_{02}} )$$}
   which, by the relation~\ref{eq:thetadet}, maps the global
section $\alpha_M^*\theta_{F_\eta}$ to $\det\beta$.

Our task now consists of relating the determinant of $\beta$ with
that of $S_M^h(x,y)$, since $\det S_M^h(x,y)=E(x,y)^{\bar r}\det
S_M(x,y)$. If fact it will be shown that the morphism $S_M^h(x,y)$
factorizes as $\beta\circ\phi^{-1}$ where $\phi$ is a morphism
whose determinant equals $\theta_{F_\eta}(M)$.

Let us begin with the morphism $\phi$. Analogous arguments as
previously applied to the exact sequence:
$$0\to \o\to \o(\Delta_{01})\to
\o(\Delta_{01})\vert_{\Delta_{01}}\to 0$$
show that there is an isomorphism:
$$ q_*(p^*M'\otimes \o(\Delta_{01}))\,\overset{\phi}\simeq\,
q_*\big((p^*M'\otimes \o(\Delta_{01})\vert_{\Delta_{01}}\big)$$ If
$\alpha_0:C\times C\to\U(r,d)$ is the morphism that sends $(x,y)$
to $M$, it then follows that the isomorphism: {\small
$$\alpha_0^*\o(\Theta_{[F_\eta]})\,\simeq\, \wedge\big(
q_*(p^*M'\otimes \o(\Delta_{01}))\big)^* \otimes \wedge
q_*(p^*M'\otimes \o(\Delta_{01})\vert_{\Delta_{01}} )\,\simeq \,
\o$$} maps $\alpha_0^*(\theta_{F_\eta})=\theta_{F_\eta}(M)$ to
$\det\phi$.

In order to write down the factorization of $S^h_M(x,y)$ we need
the following identifications: {\small
$$\begin{gathered} q_*\big(p^* M'\otimes
\o(\Delta_{01})\big)\vert_{\Delta_{02}}\,\simeq\, \iota_2^*(p^*
M'\otimes \o(\Delta_{01}))\,\simeq\, p_2^*M'\otimes
\o(\Delta) \\
   q_* \big(p^* M'\otimes
\o(\Delta_{01})\big)\vert_{\Delta_{01}}\,\simeq\, \iota_1^*(p^*
M'\otimes \o(\Delta_{01})) \, \simeq \,   p_1^*
(M'\otimes\omega_C^*)
\end{gathered}
$$}
where $\iota_j$ ($j=1,2$) is the embedding $C\times C\simeq
\Delta_{0j}\subset C\times C\times C$ and  $p_j$ is the projection
from $C\times C$ onto its $j$-th factor ($j=1,2$).

These identifications  shows that there is a natural map of
bundles on $C\times C$: {\small
$$\begin{gathered}
     p_1^*
(M'\otimes\omega_C^*) \,\simeq\, q_* \big((p^* M'\otimes
\o(\Delta_{01}))\vert_{\Delta_{01}}\big)
\,\overset{\phi^{-1}}\simeq\, q_* \big(p^* M'\otimes
\o(\Delta_{01})\big)\,\overset{\beta}\to \\
\to\,   q_* \big((p^* M'\otimes
\o(\Delta_{01}))\vert_{\Delta_{02}}\big)\,\simeq\, p_2^*M'\otimes
\o(\Delta)
\end{gathered}
$$}

If we check that this map coincides with
$S_M^h(x,y)=E(x,y)S_M(x,y)$, the lemma is proved.  To this goal it
is enough to verify that the restriction of $\beta\circ\phi^{-1}$
to the diagonal is the identity map and this fact follows from a
straightforward calculation.
\end{pf}

\begin{pf}[{\bf of Theorem~\ref{thm:add-szego}}]
Firstly, observe that Lemma~\ref{lem:sz2} implies that both sides
of the equality are holomorphic global sections of the same line
bundle.

Similar arguments to those of the proof of
Theorem~\ref{thm:add-theta} allows us to reduce the proof to check
that the statement holds true on the fiber $\pi^{-1}(z)\simeq C$
where $\pi:C^{2m}\to C^{2m-1}$  is the projection that forgets
$y_1$ and $z$ is a point $(x_1,x_2,y_2,\dots, x_m,y_m)\in
C^{2m-1}$ such that $x_i\neq x_j$ for all $i\neq j$ and $y_i=x_i$.

Now, note that the claim restricted to the fibre $\pi^{-1}(z)$ is
precisely the statement of Lemma~\ref{lem:sz3}, which has been
already proved.
\end{pf}

\begin{cor}
Under the same hypothesis of the previous theorem, one has that:
{\small
$$\begin{aligned} \operatorname{det} \left( \frac{\theta_{F_\eta
}(M(x_i-y_j ))}{\theta_{F_\eta }(M)E(x_i,y_j)^{\bar r} }\right)
\,=\, \operatorname{det} S_M(x,y)
\end{aligned}$$}
\end{cor}

\begin{pf}
It follows from Theorem~\ref{thm:add-theta} and
Theorem~\ref{thm:add-szego}.
\end{pf}

\section{Relation with the Multicomponent KP Hierarchy}

In this section, it will addressed the relation between some
properties of  non-abelian theta functions with those of
$\tau$-functions of the multicomponent KP hierarchy. The
importance of the theorem below comes from the consequences of its
infinitesimal version (Lemma~2.7 of \cite{Pl}), which eventually
leads to the bilinear identity in the framework of the
multicomponent KP hierarchy. Moreover, it generalizes
Proposition~2.16 of \cite{Fay} for higher rank.

\begin{thm}\label{thm:multicomponente}
Let $M$ be a rational point of $\U(r,r(g-1+m))$ ($m$ being a
positive integer) such that $h^1(C,M)=0$.

Then, the following identity on $C^m$ holds:
$$\theta_r(M(-\sum_{i=1}^m
y_i))\cdot\prod_{i<j}E(y_i,y_j)^r \,=\,
\lambda\cdot\det(s_i(y_j))$$ where $\lambda\in {\mathbb C}^*$,
$\{s_i=(s_i^1,\dots,s_i^r)\,\vert\, i=1,\dots,m r\}$ is a basis of
$H^0(C,M)$ and the matrix $(s_i(y_j))$ is: {\small$$ \pmatrix
s^1_1(y_1) & \ldots & s^r_1(y_1)  & \ldots & s^1_1(y_m)
& \ldots & s^r_1(y_m) \\
\vdots & & & & & & \vdots \\
            s^{1}_{m r}(y_1) & \ldots & s^r_{m r}(y_1)  & \ldots &
s^{1}_{m r}(y_m) & \ldots & s^r_{m r}(y_m)
\endpmatrix$$}
\end{thm}

\begin{pf}
We begin with the $m=2$ case where the idea of the proof will be
clear. For this case we will work with bundles on $C\times \times
C\times C$ and will use again the notations introduced in the
proof of Lemma~\ref{lem:sz3}.

Recall that the sheaf $\o_{\Delta_{01}+\Delta_{02}}$ is  the
kernel of the difference map
$\o_{\Delta_{01}}\oplus\o_{\Delta_{02}}\to
\o_{\Delta_{01}\cap\Delta_{02}}$ and, therefore, we have the exact
sequence:
$$0\to \o_{\Delta_{01}+\Delta_{02}}\to
\o_{\Delta_{01}}\oplus\o_{\Delta_{02}}\to
\o_{\Delta_{01}\cap\Delta_{02}}\to 0$$

  From the following exact  sequence:
$$0\to \o(-\Delta_{01}-\Delta_{02})\to \o
\to \o_{\Delta_{01}+\Delta_{02}}\to 0$$ one deduces the exactness
of: {\small
$$\begin{aligned}
   0 \to q_*(p^*M(-\Delta_{01}-\Delta_{02}))\to &
q_*(p^*M)\overset{\alpha}\to \\
\to q_*(p^*M\otimes \o_{\Delta_{01}+\Delta_{02}})\to  & R^1
q_*(p^*M(-\Delta_{01}-\Delta_{02}))\to 0 \end{aligned}$$} which,
by the relation~\ref{eq:thetadet}, implies that:
$$\theta_r(M(-y_1-y_2))=\det(\alpha)\qquad\forall y_1,y_2\in C$$

The statement is thus reduced to compute $\det(\alpha)$ in an alternative
way.

Consider the following commutative diagram: $$\xymatrix{ q_*(p^*M)
\ar[r]^(.4){\alpha} \ar@{=}[d] &
     q_*(p^*M\otimes \o_{\Delta_{01}+\Delta_{02}}) \ar[d]^{\varphi} \\
H^0(M)\otimes \o_{C\times C} \ar[r]^(.4){\text{\it ev}} &
q_*(p^*M\otimes (\o_{\Delta_{01}}\oplus\o_{\Delta_{02}}))}$$ where
$\varphi$ is the morphism induced by
$\o_{\Delta_{01}+\Delta_{02}}\to
\o_{\Delta_{01}}\oplus\o_{\Delta_{02}}$ and $\text{\it ev}$ is the
evaluation map, that is, at the point $(y_1,y_2)$ is:
$$\begin{aligned}
H^0(M)\,&\longrightarrow\, M_{y_1}\oplus M_{y_2} \\
s\,&\longmapsto (s(y_1),s(y_2))
\end{aligned}$$

The diagram shows that:
$$\det(\varphi)\cdot\det(\alpha)\,=\, \det(\text{\it ev})$$
and therefore:
$$\det(\varphi)\cdot\theta_r(M(-y_1-y_2))\,=\, \lambda'\det(s_i(y_j))$$
where $\lambda'$ is a constant that depends on the choice of the
basis and on the above isomorphisms of line bundles and it will
eventually give the constant of the statement.

Since $q: \Delta_{01}+\Delta_{02}\to C$ is finite of degree $2$,
   $R^1 q_*(\o_{\Delta_{01}+\Delta_{02}})=0$. It thus
follows the exactness of:
$$0\to q_*(\o_{\Delta_{01}+\Delta_{02}})\overset{\varphi_0}\to
q_*(\o_{\Delta_{01}}\oplus\o_{\Delta_{02}}) \to q_*(
\o_{\Delta_{01}\cap \Delta_{02}})\to 0$$ Now, we will show that
$\det(\varphi)=\det(\varphi_0)^r$ and that
$\det(\varphi_0)=E(y_1,y_2)$.

Let us begin computing $\det(\varphi_0)$. From the theory of
determinants (\cite{KM}) one has the following isomorphism:
$$\Det\big(
q_*( \o_{\Delta_{01}+\Delta_{02}})\to q_*
(\o_{\Delta_{01}}\oplus\o_{\Delta_{02}})\big) \,\simeq\, \Det(q_*(
\o_{\Delta_{01}\cap \Delta_{02}}))$$

Since the above bundles live on $C\times C$ let us rewrite them as
follows. From the diagram: {\small $$\xymatrix{ \Delta_{01}\cap
\Delta_{02} \ar[d] \ar[r] & C\times C\times C \ar[d]^q \\ \Delta
\ar[r] & C\times C }$$} ($\Delta\subset C\times C$ being the
diagonal) one obtains:
$$\Det(q_*( \o_{\Delta_{01}\cap \Delta_{02}}))
\,\simeq\, \Det\big(\o_{\Delta}\big) \,\simeq\,
\Det\big(\o(-\Delta)\to \o\big) \,=\, \o(-\Delta)$$ where the
second isomorphism follows from the exactness of the following
sequence on $C\times C$:
$$0\to \o(-\Delta)\to \o \to
\o_{\Delta}\to 0$$

These calculations imply that:
$$\det(\varphi_0)=E(y_1,y_2)\,\in \, H^0(C\times C,\o(\Delta))$$

On the other hand, $\det(\varphi)$  may be computed similarly and
we obtain:
$$\begin{aligned}
\Det\big( q_*(p^*M\otimes  & \o_{\Delta_{01}+\Delta_{02}})\to q_*
(p^*M\otimes
(\o_{\Delta_{01}}\oplus\o_{\Delta_{02}}))\big)\,\simeq
\\ &
\simeq\, \Det\big(M\otimes \o_{\Delta}\big) \,\simeq\,
   \o(-r\Delta) \end{aligned}$$

   Now, it follows that
   $\det(\varphi)=\det(\varphi_0)^r=E(y_1,y_2)^r$.
The $m=2$ case is proved.

For arbitrary $m$ we proceed similarly but replacing the morphism
$\varphi$ by: {\small $$\begin{aligned} q_*(p^*M\otimes
\o_{\sum_{i=1}^m\Delta_{0i}})\,\to &\, q_*(p^*M\otimes
(\o_{\sum_{i=1}^{m-1}\Delta_{0i}}\oplus \o_{\Delta_{0m}}))
\,\to\,\dots \\
&\dots \,\to\, q_*(p^*M\otimes
(\oplus_{i=1}^m\o_{\Delta_{0i}}))\end{aligned}$$} which has
determinant $\prod_{i<j}E(y_i,y_j)^r$.
\end{pf}

\section{Cyclic coverings}

Let $\gamma:\wtilde C\to C$ be a cyclic covering of degree $n$
between two irreducible smooth projective cuves given by an
automorphism $\sigma$ of $\wtilde C$ such that
$\sigma^n=\operatorname{Id}$, that is, $\wtilde C/<\sigma>= C$.

In this section we will study the relationship between  the
polarizations of moduli spaces of vector bundles on $\wtilde C$
and $C$. This question is related to the twist structures of $b-c$
systems (\cite{Raina-bc}) and has been addressed in \cite{Ball}.
The rank $1$ case is to be found in \cite{Fay}.

Let us introduce some notation. Let $\Delta$ (resp. $\wtilde
\Delta$) denote the diagonal of $C\times C$ (resp. $\wtilde
C\times\wtilde C$). Let $\wtilde D_{ij}$ be the inverse image of
the diagonal by the morphism $\sigma^i\times\sigma^j:\wtilde
C\times \wtilde C\to\wtilde C\times \wtilde C$. Let $E$ (resp.
$\wtilde E$) be the prime form of $C$ (resp. $\wtilde C$).
Finally, let $R^\gamma=\sum_{\tilde x\in\wtilde C}(n_{\tilde
x}-1)\wtilde x$ be the ramification divisor of $\gamma$, where
$n_{\tilde x}$ is the ramification index at $\wtilde x$.


Let us begin with some computations for the ideal sheaf of the
diagonal.

\begin{lem}\label{lem:diagonal}
Let $\gamma_1$ be $\gamma\times\gamma$ and  $R^1$ be
$(R^\gamma\times\wtilde C)\cup (\wtilde C\times R^\gamma)$.

Then, there is an exact sequence on $\wtilde C\times\wtilde C$:
$$0\to \gamma_1^*\o(\Delta)\to \o(\sum_{j}\wtilde D_{0j})\to
           \o_{R^1}\to 0$$
and a canonical isomorphism of line bundles:
$$ \gamma_1^*\o(\Delta) \,\simeq \,
\o(\sum_{j}\wtilde D_{0j})\otimes \wtilde
p_1^*\o(-R^\gamma)\otimes \wtilde p_2^*\o(-R^\gamma)$$ where
$\wtilde p_i:\wtilde C\times\wtilde C\to\wtilde C$ are the natural
projections.
\end{lem}

\begin{pf}
Since there is an inclusion $\gamma_1^{-1}(\Delta)\subseteq
\sum_j \wtilde D_{0j}$, it follows the  exact
sequence:
\beq
0\to \gamma_1^*\o(\Delta)\to \o(\sum_{j}\wtilde D_{0j})\to
           \o_T\to 0
\label{eq:suc1}\end{equation}

Let us compute  $\o_T$. If $S=\supp(R^\gamma)$ is the support of $R^\gamma$
and $U$ is the open subscheme $(\wtilde C-S)\times (\wtilde C-S)$, one
checks easily that
$\gamma_1^{-1}(\Delta)\vert_U = \sum_j\wtilde D_{0j}\vert_U$ and $T$ is
therefore contained in:
           $$\wtilde C\times\wtilde C-U=\underset{\tilde x\in
R^\gamma}\bigcup\big(
\{\wtilde x\}\times\wtilde C\cup \wtilde C\times \{\wtilde x\}\big)$$

By symmetry, it is enough to show that the length of $T$ at
$\{\wtilde x\}\times\wtilde C$ ($\wtilde x\in S$) is $n_{\tilde
x}-1$. Recalling the exact sequence \ref{eq:suc1}, one observes
that this can be done by comparing the zero divisors of
$\gamma_1^*E$ and $\prod_j
(\operatorname{Id}\times\sigma^j)^*\wtilde E$ as global sections
of $\o(\sum_{j}\wtilde D_{0j})$. One checks now that if $(\wtilde
x,\wtilde y) \in S\times \wtilde C$, then $(\wtilde x,\wtilde y)$
           is a simple zero of $\gamma_1^*E$ and a zero of order
$n_{\tilde x}$ of
$\prod_j (\operatorname{Id}\times\sigma^j)^*\wtilde E$.

For the second claim, it suffices to take determinants in the
exact sequence of the first claim.
\end{pf}

\begin{lem}\label{lem:wedge}
Let $\wtilde M$ be a vector bundle on $\wtilde C$ of rank $\tilde
r$.

Then, there is an exact sequence:
$$0\to\gamma^*(\gamma_*\wtilde M)\to
\oplus_{k=0}^{n-1} (\sigma^k)^*\wtilde M\to \big(\o_{\frac12 n
R^\gamma}\big)^{\oplus\tilde r}\to 0$$ and a canonical isomorphism:
$$\wedge \gamma^*(\gamma_*\wtilde M)\,\simeq\, \otimes_{k=0}^{n-1}\wedge
(\sigma^k)^*\wtilde M\otimes \o(-\frac12{\tilde  r}n R^\gamma)$$

Finally, if $\wtilde M$ has degree $\wtilde d$, then:
$$\deg\gamma_*\wtilde M\,=\, {\wtilde  d}-{\wtilde  r}(\wtilde
g-1-n(g-1))$$ where $\wtilde g$ (resp. $g$) is the genus of $\wtilde C$
(resp. $C$).
\end{lem}

\begin{pf}
The adjunction formula gives a morphism $\gamma^*(\gamma_*\wtilde M)\to
\wtilde M$ and,
since $\gamma^*(\gamma_*\wtilde M)$ is invariant under $\sigma$, there is
also a
morphism to $(\sigma^k)^*\wtilde M$ for  $0\leq k< n$; that is:
$$\gamma^*(\gamma_*\wtilde M)\to \oplus_{k=0}^{n-1} (\sigma^k)^*\wtilde M$$

Since this is a morphism between two locally free sheaves of the same rank
which is an isomorphism at the stalk of any point $\wtilde x\in\wtilde
C-R^\gamma$, it follows that there is an exact sequence:
\beq
0\to\gamma^*(\gamma_*\wtilde M)\to
\oplus_{k=0}^{n-1} (\sigma^k)^*\wtilde M\to \o_T\to 0
\label{eq:suc2}\end{equation}
           where $\supp(T)\subseteq\supp R^\gamma$.

Note that the computation of $T$ is a local problem, so it can be assumed
$\wtilde M$ to be $\o_C^{\oplus\tilde r}$. Furthermore, observe that
$\o_T=\o_{T'}^{\oplus\tilde r}$ where:
\beq
0\to\gamma^*(\gamma_*\o_{\wtilde C})\to
\oplus_{k=0}^{n-1} \o_{\wtilde C}\to \o_{T'}\to 0
\label{eq:trivial}\end{equation}

For the case $\wtilde M=\o_{\wtilde C}$ some results on
cyclic coverings are known. From Theorem~3.2 of \cite{esteban} we learn
that the covering
$\gamma:\wtilde C\to C$ is defined by a line bundle $L$ on $C$ and a
divisor $D=\sum a_iq_i$ on $C$ where: $1\leq a_i< n$, $L^{\otimes n}\simeq
\o_C(D)$, and $q_i$ is a branch point of $\gamma$. Furthermore,
all the points on the fibre of a $q_i$ have the same multiplicity, say
$m_i$, and $s_i:=\frac{n}{m_i}=\text{g.c.d}(a_i,n)$ is  the
number of distinct points in $\gamma^{-1}(q_i)$.

Moreover, if $[a]_n$ denotes the remainder of $a$ divided by $n$ and $D_k$
is
$\sum_i [k a_i]_n q_i$, it then holds that the
coefficients of $\gamma^{-1}(D_k)$ are multiple of $n$ (\S2 of
\cite{esteban}) and that:
$$ \gamma^*(\gamma_*\o_{\wtilde C})\,\overset{\sim}\longrightarrow \,
\oplus_{k=0}^{n-1} \o_{\wtilde C}(-\frac{1}{n}\gamma^{-1}(D_k))$$

Now, one checks that the  morphism \ref{eq:trivial} is given by
the divisors $-\frac{1}{n}\gamma^{-1}(D_k)$, in particular,
$\supp(T')\subseteq \cup_k\supp \gamma^{-1}(D_k)=\supp R^\gamma$.

It only remains to compute the length of $T'$ at a ramification
point. Let $p_i\in\gamma^{-1}(q_i)$ be given.

The length of the cokernel of the sequence \ref{eq:trivial} at
$p_i$ is given by: {\small
$$\begin{aligned}\sum_{k=1}^{n-1}[k a_i]_n\frac{ m_i}{n} \,=\, &
\frac{m_i s_i}{n}\sum_{k=1}^{m_i-1}[k a_i]_n \,=\, s_i
\sum_{k'=1}^{m_i-1}[k']_n \,=\,\\  \,=\, & s_i
\frac{m_i(m_i-1)}{2}\,=\, \frac{n(m_i-1)}{2}\end{aligned}$$} Thus,
$T'=\frac12{n} R^\gamma$ and the conclusion follows. Observe that
the coefficients of $n R^\gamma$ are even.
\end{pf}


Fix a line bundle $L_\gamma$ on $\wtilde C$ satisfying:
$$\cases
L_\gamma\,:=\,\o_{\wtilde C}(\frac12 (n R^\gamma)-m R^\gamma) &
\text{if }n=2m+1 \\
L_\gamma\text{ such that } L_\gamma^{\otimes 2}\simeq \o_{\wtilde
C}(R^\gamma) &
\text{if }n=2m
\endcases$$
Then, the following two conditions hold:
$$\begin{aligned}
L_\gamma^{\otimes 2}\,\simeq&\, \o_{\wtilde C}(R^\gamma)\\
L_\gamma^{\otimes n}\,\simeq&\, \o_{\wtilde C}(\frac12 (n
R^\gamma))\end{aligned}$$ and $L_\gamma$ has degree $(\wtilde g
-1)-n(g-1)$.

We also fix theta characteristics $\eta$ on $C$ and $\wtilde \eta$
on $\wtilde C$ where $\wtilde \eta$ is defined by $\o_{\wtilde
C}(\wtilde\eta):=\gamma^*\o_C(\eta)\otimes L_\gamma$.

Let $\wtilde d:= \wtilde r(\wtilde g-1-n(g-1))$. Since
$p.g.c.d.(\wtilde r,\wtilde d)=\wtilde r$, the line bundle
$F=L_{\gamma}^*$ may be used to define a polarization
$\wtilde\Theta_{[F_{\tilde\eta}]}$ in $\U_{\wtilde C}(\wtilde
r,\wtilde d)$.

Note that the theta characteristic $\eta$ also defines a
polarization $\Theta_{[\eta]}$ on the moduli space $\U_C(r,0)$.

Assume, we are given a vector bundle $\wtilde M\in \U_{\wtilde
C}(\wtilde r,\wtilde d)$ whose direct image is a semistable vector
bundle on $C$. Then,  Lemma~\ref{lem:wedge} implies that
$\gamma_*\wtilde M\in\U(r,0)$ where $r:=n\cdot\wtilde r$. Further,
we have the morphisms:
$$\begin{aligned}
\wtilde \alpha_{\wtilde M}\colon \wtilde C^{2nm} \,&
\longrightarrow\,
\U_{\wtilde C}(\wtilde r,\wtilde d) \\
(\tilde x_1,\tilde y_1,\dots,\tilde x_{nm},\tilde
y_{nm})&\longmapsto \wtilde M(\sum_i \tilde x_i-\tilde y_i)
\end{aligned}$$
and
$$\begin{aligned}
   \alpha_{\gamma_*\wtilde M}\colon  C^{2m} \,& \longrightarrow\,
\U_{ C}( r,0) \\
(x_1,y_1,\dots,x_{m},y_{m})&\longmapsto \gamma_* \wtilde
M\otimes\o_C(\sum_i x_i-y_i)
\end{aligned}$$

   In order to study the relation
of the corresponding non-abelian theta functions, we consider the
following diagram:

$$\xymatrix{
& (\prod^n (\wtilde C\times\wtilde C))^{m}=\wtilde C^{2nm}
\ar[r]^-{\tilde\alpha_{\wtilde M}} & \U_{\wtilde C}(\wtilde r,\wtilde d)\\
(\wtilde C\times\wtilde C)^{m} \ar[ru]^-{\rho_m} \ar[rd]_-{\gamma_m} \\
& (C\times C)^{m}=C^{2m} \ar[r]^-{\alpha_{\gamma_*\wtilde M}} &
\U_{C}(r,0) }$$ where $\gamma_m$ denotes the map $\wtilde
C^{2m}\to C^{2m}$ given by $\gamma$ on each component, and
$\rho_m$ is the embedding induced by the morphism:
$$\begin{aligned}
\rho_1:\wtilde C\times\wtilde C\,& \longrightarrow\, \,
\prod^n(\wtilde C\times\wtilde C) \\
(\wtilde x,\wtilde y)\,& \longmapsto (\wtilde x,\wtilde
y,\sigma(\wtilde x),\sigma(\wtilde y), \dots,\sigma^{n-1}(\wtilde
x),,\sigma^{n-1}(\wtilde y))
\end{aligned}
$$

Let us denote by $p_i$ (resp. $\wtilde p_i$) the projection of
$C^{2m}$ (resp. $\wtilde C^{2m}$) onto its $i$-th factor.

The following theorem gives the relation between the pullbacks of
the polarizations  by the above diagram.

\begin{thm}\label{thm:cyclic-direct}
There is an isomorphism of line bundles on $\wtilde C^{2m}$:
{\small $$(\alpha_{\gamma_*\wtilde M} \circ
\gamma_m)^*\o(\Theta_{[\eta]}) \,\simeq \,
            \rho_m^*\left(
(\wtilde\alpha_{\wtilde
M})^*\o(\wtilde\Theta_{[F_{\tilde\eta}]})\right) \otimes
\big(\underset{i\text{ all}}\otimes{\wtilde p_i^*}
L_\gamma^*\big)^{\otimes 3r}
$$}
where $F=L_\gamma^*$ and $r=\wtilde r n$.
\end{thm}

\begin{pf}
            The statement follows from
the comparation of the pullbacks  $(\alpha_{\gamma_*\wtilde M}
\circ \gamma_m)^*\o(\Theta_{[\eta]})$ and $(\wtilde\alpha_{\wtilde
M}\circ \rho_m)^*\o(\wtilde\Theta_{[F_{\tilde\eta}]})$, which will
be done with the help of Theorem~\ref{thm:main-iso} and
Lemmas~\ref{lem:diagonal} and \ref{lem:wedge}.

To begin with, we compute the pull-back by $\gamma_m$ of
$(\alpha_{\gamma_*\wtilde M} )^*\o(\Theta_{[\eta]})$. Note
that $L_\gamma^*(-\wtilde \eta)$ is invariant by $\sigma$. Therefore, by
Lemma
\ref{lem:wedge} and the properties of $L_\gamma$, one has that:
   {\small $$\begin{aligned}
   \gamma&_m^* \left(
\underset{i\text{ odd}}\otimes p_i^* \wedge \left( \gamma_*\wtilde
M\otimes\o(-\eta)\right)^*\right)
\,\simeq \,
\underset{i\text{
odd}}\otimes\wtilde p_i^*\left(\gamma^* \wedge
\left( \gamma_*\wtilde M\otimes\o(-\eta)\right)^*\right)
\,\simeq \\
   & \simeq\,\underset{i\text{ odd}}\otimes\wtilde p_i^*
\wedge \left( \gamma^*\gamma_*\wtilde
M\otimes\gamma^*\o(-\eta)\right)^*
   \,\simeq \,
   \underset{i\text{ odd}}\otimes\wtilde p_i^*\left( \wedge
\gamma^*\gamma_*\wtilde M\otimes L_\gamma(-\wtilde \eta)^{\otimes
r}\right)^* \,\simeq \\ & \simeq\, \underset{i\text{
odd}}\otimes\wtilde p_i^*\left( \big( \otimes_{j=0}^{n-1}\wedge
(\sigma^j)^* \wtilde M\big) \otimes \o(-\frac12 r R^\gamma)
\otimes L_\gamma(-\wtilde \eta)^{\otimes r}\right)^*
\end{aligned}$$}
Recalling that $(L_\gamma^*)^{\otimes r}\simeq\o(-\frac12 r
R^\gamma)$ and that $L_\gamma^*(-\wtilde \eta)$ is invariant under
$\sigma$, the above expression is isomorphic to:
  {\small $$\begin{aligned}
\underset{i\text{ odd}}\otimes & \wtilde p_i^*\left( \big(
\otimes_{j=0}^{n-1}\wedge (\sigma^j)^* \wtilde M\big) \otimes
L_\gamma^*(-\wtilde \eta)^{\otimes r}\otimes L_\gamma^{\otimes r}
\right)^* \,\simeq \\ & \simeq\,
   \underset{i\text{ odd}}\otimes\wtilde p_i^*\left( \left(
\otimes_{j=0}^{n-1}\wedge (\sigma^j)^* (\wtilde M \otimes
L_\gamma^*(-\wtilde \eta))\right) \otimes L_\gamma^{\otimes r}
\right)^* \,\simeq \\ & \simeq\, \rho_m^*\left( \underset{i\text{
odd}}\otimes{\overline p_i^*} \wedge \big(\wtilde M\otimes
L_\gamma^* (-\wtilde \eta)\big)^*\right) \otimes \big(
\underset{i\text{ odd}}\otimes  \wtilde p_i^*(L_\gamma^*)^{\otimes
r}\big)
\end{aligned}$$} where $\overline p_i$ are the natural projections of
$\wtilde C^{2nm}$.

Similarly, the pullback of: {\small
$$\underset{i\text{ even}}\otimes p_i^*\wedge
\left(\gamma_*\wtilde M\otimes\o(\eta)\right)$$}by $\gamma_m^*$
is:
$$\rho_m^*\left( \underset{i\text{
even}}\otimes{\overline p_i^*}
\wedge \big(\wtilde M\otimes L_\gamma^* (\wtilde
\eta)\big)\right) \otimes \big( \underset{i\text{ even}}\otimes\wtilde
p_i^*(L_\gamma^*)^{\otimes r}\big)$$

Note that $n\sum_i\wtilde D_{0i}=\sum_{i,j}\wtilde D_{ij}$ and
$\rho_1^{-1}(\wtilde \Delta_{ij})=\wtilde D_{ij}$. Then, from
Lemma \ref{lem:diagonal}, it follows an isomorphism on $\wtilde
C\times \wtilde C$:
$$
\begin{aligned}
\gamma_1^*\o(n\Delta)\,& \simeq\, \rho_1^*\o(\sum_{i,j}^n
\wtilde\Delta_{ij})\otimes \wtilde p_1^*\o(-nR^\gamma)\otimes
\wtilde p_2^*\o(-nR^\gamma)
\\ & \simeq\,
  \rho_1^*\o(\sum_{i,j}^n
\wtilde\Delta_{ij})\otimes \wtilde p_1^*(L_\gamma^*)^{\otimes 2n}
\otimes \wtilde p_2^*(L_\gamma^*)^{\otimes 2n}\end{aligned}
$$
Finally, a length but straightforward calculation shows that:
{\small
$$\begin{aligned} \gamma_m^*\o(\sum_{\Sb i+j=\text{odd}\\ i<j
\endSb}\Delta_{ij}  &
-\sum_{\Sb i+j=\text{even}\\ i<j \endSb}\Delta_{ij})^{\otimes r} \,\simeq\\
\simeq\, & \rho_m^*\o\Big(\sum_{\Sb i+j=\text{odd}\\ i<j
\endSb}\wtilde\Delta_{ij}-
\sum_{\Sb i+j=\text{even}\\ i<j
\endSb}\wtilde\Delta_{ij}\Big)^{\otimes \wtilde r}
\otimes\big(\underset{i \text{ all}}\otimes \wtilde
p_i^*(L_\gamma^*)^{\otimes 2r}\big)
\end{aligned}
$$}

Comparing these results with the expression of
$\rho_m^*\big(\wtilde\alpha_{\wtilde
M}^*\o(\wtilde\Theta_{[F_{\tilde\eta}]})\big)$ given by
Theorem~\ref{thm:main-iso}, the statement follows.
\end{pf}

Observe that  $L_\gamma=\o_{\wtilde C}$ when $\gamma$ is
non-ramified. Then, in this situation, a consequence of the above
theorem is the following identity between global sections of the
line bundles in the previous statement:

\begin{thm}\label{thm:direct} Let $\gamma$ be non-ramified and $\wtilde
M\in \U(\wtilde r, \wtilde d)-\wtilde \Theta_{\tilde \eta}$ such
that $\gamma_*\wtilde M\in \U(r,0)$. Then, for $(\wtilde x_1,
\wtilde y_1, \dots ,\wtilde x_m,\wtilde y_m)\in \wtilde C^{2m}$,
it holds that: {\small
$$
\frac{\theta_{\eta}\left(\gamma_*\wtilde M
\otimes\o(Z)\right)}{\theta_{\eta}(\gamma_*\wtilde M)} \,=\,
\frac{\wtilde\theta_{\tilde \eta}\left(\wtilde M\otimes
\gamma^*\o(Z)\right)} {\wtilde\theta_{\tilde \eta}(\wtilde M)}
$$}
where $Z$ is the divisor $\sum_{i=1}^m(\gamma (\wtilde x_i)-\gamma
(\wtilde y_i))$ on $C$.
\end{thm}

\begin{pf}
First of all, observe that: {\small $$\wtilde\theta_{\tilde
\eta}\left(\wtilde M\otimes
\gamma^*\o(Z)\right) \,=\,
\wtilde\theta_{\tilde \eta}\left(\wtilde \alpha_{\tilde
M}(\rho_m(\wtilde x_1,\wtilde y_1,\dots, \wtilde x_m,\wtilde
y_m))\right)
$$}
because $\gamma^{-1}(\gamma(\wtilde x_i))=\rho_1(\wtilde x_i)$.
Then, the r.h.s. of the formula is a holomorphic global section of
$\rho_m^*(\wtilde \alpha_{\tilde M}^*\o(\wtilde\Theta_{[\tilde
\eta]} ))$. On the other hand, the l.h.s. is a holomorphic global
section of $\gamma_m^*(\alpha_{\gamma_*\tilde
M}^*\o(\Theta_{[\eta]}))$. Hence, by
Theorem~\ref{thm:cyclic-direct}, both sides are global section of
isomorphic line bundles on $\wtilde C^{2m}$.

Similar arguments to those of the proof of
Theorem~\ref{thm:add-theta} reduce the proof to check that the
statement holds true when restricted to a fibre $\pi^{-1}(z)$
where $\pi:\wtilde C^{2m}\to \wtilde C^{2m-1}$ is the projection
that forgets $\wtilde y_1$ and $z$ is a point $(\wtilde
x_1,\wtilde x_2,\wtilde y_2,\dots,\wtilde x_m,\wtilde y_m)\in
\wtilde C^{2m-1}$ such that $x_i\neq x_j$ for all $i\neq j$ and
$y_i=x_i$. For the sake of notation, we define $x$ to be
$\gamma(\wtilde x_1)$.

Let denote by $\wtilde p$ and $\wtilde q$ the projections of
$\wtilde C\times \wtilde C$ onto its first and second factors,
respectively. Consider the  bundle:
$$\wtilde\M:=\wtilde
p^*(\wtilde M(\gamma^{-1}(x)+\wtilde\eta)\otimes \o(-\wtilde D)$$
on $\wtilde C\times\wtilde C$, where $\wtilde D:=\sum_j\wtilde
D_{0j}$. Using the sequence defined by the effective divisor
$\wtilde D$, we obtain the following exact sequence: {\small
\beq\label{eq:mtilde} 0\to \wtilde \M\to \wtilde p^*(\wtilde
M(\gamma^{-1}(x)+\wtilde\eta))\to \wtilde p^*(\wtilde
M(\gamma^{-1}(x)+\wtilde\eta))\otimes\o_{\wtilde D}\to 0
\end{equation}} Since $R^1 q_* (\wtilde p^*(\wtilde
M(\gamma^{-1}(x)+\wtilde\eta)))=0$, the restriction of the r.h.s.
to $\pi^{-1}(z)$ is given by the determinant of the morphism:
{\small $$ \wtilde q_*(\wtilde p^*(\wtilde
M(\gamma^{-1}(x)+\wtilde\eta)))\,\to\,  \wtilde q_*(\wtilde
p^*(\wtilde M(\gamma^{-1}(x)+\wtilde\eta))\otimes\o_{\wtilde
D})$$}
   induced by
the latter sequence.

We now compute the restriction of the l.h.s. in a similar way. Let
denote by $p$ and $q$ the projections of $C\times \wtilde C$ onto
its first and second factors, respectively. Let $\M$ be the bundle
$p^*(\gamma_*\wtilde M\otimes\o(x+\eta))\otimes\o(-\Gamma)$ on
$C\times \wtilde C$, where $\Gamma$ is the graph of the map
$\gamma$. The exact sequence associated to the divisor $\Gamma$
implies the exactness of the sequence: {\small \beq\label{eq:m}
0\to \M \to
   p^*(\gamma_*\wtilde M\otimes\o(x+\eta)) \to p^*(\gamma_*\wtilde
   M\otimes\o(x+\eta))\otimes\o_{\Gamma}\to 0
   \end{equation}}
Being
$R^1q_*(p^*(\gamma_*\wtilde M\otimes\o(x+\eta)))=0$, it follows
that the restriction of the l.h.s. to $\pi^{-1}(z)$ is the
determinant of the induced morphism: {\small $$
q_*(p^*(\gamma_*\wtilde M\otimes\o(x+\eta)))\,\to \,
q_*(p^*(\gamma_*\wtilde M\otimes\o(x+\eta))\otimes\o_{\Gamma})$$}

Bearing in mind the commutativity of the diagram:
$$\xymatrix{\wtilde C\times\wtilde C \ar[dr]_{\wtilde q}
\ar[rr]^{\gamma\times \operatorname{Id}} & & C\times \wtilde C
\ar[dl]^q \\
& \wtilde C}$$ it will suffice to conclude to show that the direct
image by $\gamma\times \operatorname{Id}$ of the
sequence~\ref{eq:mtilde} is the sequence~\ref{eq:m}.

The direct image of the sequence~\ref{eq:mtilde} by $\gamma\times
\operatorname{Id}$ is: {\small $$\begin{aligned}
  0\to  (\gamma\times
\operatorname{Id})_*\Big(\wtilde p^*(\wtilde M(\gamma^{-1}(x) &
+\wtilde\eta))\otimes\o(-{\wtilde D})\Big)\to\\  \to (\gamma\times
\operatorname{Id})_*& \Big(\wtilde p^*(\wtilde
M(\gamma^{-1}(x)+\wtilde\eta))\Big)\to\\  &\to (\gamma\times
\operatorname{Id})_*\Big(\wtilde p^*(\wtilde
M(\gamma^{-1}(x)+\wtilde\eta))\otimes\o_{\wtilde D}\Big)\to 0
\end{aligned}$$}
because the map $\gamma\times \operatorname{Id}$ is finite.

Recalling that $(\gamma\times
\operatorname{Id})^{-1}(\Gamma)=\wtilde D$,
$\gamma^{-1}(\eta)=\wtilde\eta$ and using the projection formula
and the base change theorem for the case:
$$\xymatrix{\wtilde C\times\wtilde C \ar[d]_{\wtilde p} \ar[r]^{\gamma\times
\operatorname{Id}} & C\times\wtilde C \ar[d]^p \\ \wtilde C
\ar[r]^\gamma & C}$$
  we conclude that the latter sequence
coincides with the sequence~\ref{eq:m}.

Then, we know that both sections are equal up to a constant on
$\wtilde C^{2m}$. This constant might be evaluated on
$\pi^{-1}(z)$ by letting $\wtilde y_1=\wtilde x_1$, and it follows
that it is equal to $1$.
\end{pf}

\begin{rem}
Observe that the above theorem may be generalized for the ramified
case. This would require to know sections of $L_\gamma^{\otimes
n}=\o(\frac12 n R^\gamma)$. Besides, Lemma~\ref{lem:diagonal}
allows us to give a section of it in terms of the prime forms $E$
and $\wtilde E$ because $\frac12 n R^\gamma$ is effective.
\end{rem}


We finish with a similar study for the inverse image. Let  $M\in
\U_C(r,0)$ be a vector bundle on $C$. From Lemma~3.2.2 of
\cite{ML} it turns out that $\gamma^* M\in\U_{\wtilde C}(r,0)$.
Consider the following diagram:

$$\xymatrix{
           {\wtilde C^{2m}} \ar[d]_-{\gamma_m}
\ar[r]^-{\tilde\alpha_{\gamma^* M}} &
\U_{\wtilde C}(r,0)\\
C^{2m} \ar[r]^-{\alpha_M} & \U_{C}(r,0) }$$

\begin{thm}\label{thm:cyclic-inverse}
It holds that: {\small $$\begin{aligned} (\alpha_M &\circ
\gamma_m)^* \o(\Theta_{[\eta]})\otimes
\gamma_m^*\o\big(\sum_{i+j=\text{odd}}
\Delta_{ij}-\sum_{i+j=\text{even}}\Delta_{ij}\big)^{\otimes r}
\,\simeq \\ \simeq\, & (\tilde\alpha_{\gamma^*
M})^*\o(\wtilde\Theta_{[\tilde\eta]})\otimes
\o(\sum_{i+j=\text{odd}}\wtilde
\Delta_{ij}-\sum_{i+j=\text{even}}\wtilde \Delta_{ij})^{\otimes r}
\otimes\big(\underset{i\text{ all}}\otimes (\wtilde p_i^*
L_\gamma^*)^{\otimes r} \big)
\end{aligned}$$}
\end{thm}

\begin{pf}
The claim follows from Theorem \ref{thm:main-iso} and from the
fact that $\gamma^*\o_C(\eta ) \simeq \o_{\wtilde C}(\wtilde
\eta)\otimes L_\gamma^*$.
\end{pf}

\begin{prop}\label{prop:inverseimage}
Suppose that $\gamma$ is non-ramified and let  $M\in \U_C(r,0)$
such that $\theta_{\eta}(M)\neq 0$. Then, for $(\wtilde
x_1,\wtilde y_1,\dots,\wtilde x_m,\wtilde y_m)\in \wtilde C^{2m}$,
it holds that: {\small
$$\begin{gathered}  \frac{\wtilde\theta_{\tilde\eta}
\Big(\gamma^*M (\sum_{i=1}^m(\wtilde x_i-\wtilde
y_i))\Big)}{\wtilde\theta_{\tilde\eta}(\gamma^*M)}\cdot
\prod_{i<j}\prod_{k=1}^{n-1} \Big( \wtilde E(\wtilde
x_i,\sigma^k(\wtilde x_j))  \wtilde E(\wtilde y_i,\sigma^k(\wtilde
y_j)) \Big)^r
  \,=\\ \,\qquad \qquad \qquad\qquad =\,
\frac {\theta_\eta\Big(M(\sum_{i=1}^m(\gamma(\wtilde
x_i)-\gamma(\wtilde y_i)))\Big)} {\theta_\eta(M)} \cdot
  \prod_{i,j}
\prod_{k=1}^{n-1} \wtilde E(\wtilde x_i,\sigma^k(\wtilde y_j))^r
\end{gathered}$$}
\end{prop}

\begin{pf}
One proceed similarly as in Theorem~\ref{thm:direct}.
\end{pf}

\begin{rem}\label{rem:prop5.1defay}
It is worth pointing out that Proposition~5.1 of \cite{Fay}
follows from Theorem~\ref{thm:cyclic-inverse} when $\gamma$ is
ramified, $\deg\gamma=2$ and $r=1$.
\end{rem}





\vskip2truecm
\begin{thebibliography}{IMM}

\bibitem{Ball} {\it Ball, J.A.; Vinnikov, V.}, ``Zero-pole
interpolation for matrix meromorphic functions on a compact Riemann
surface and a matrix Fay trisecant identity'', Am. J. Math. {\bf 121}
(1999), pp. 841--888

\bibitem{BNR} {\it Beauville, A.; Narasimhan, M.S.; Ramanan, S.},
``Spectral covers and the generalized theta divisor'', J. reine angew.
Math. {\bf 398} (1989), pp. 169--179

\bibitem{DN} {\it Drezet, J.M.; Narashimhan, M.S.},
``Groupe de Picard des vari\'et\'es de modules de fibr\'es semi--stables
sur les courbes alg\'ebriques'',  Invent. Math. {\bf 97} (1989), pp. 53--94

\bibitem{Fay} {\it Fay, J.D.}, ``Theta Functions on Riemann Surfaces'',
LNM~352, Springer-Verlag (1973)

\bibitem{Fay2} {\it Fay, J.D.}, ``The non-abelian Szeg\"o kernel and
theta divisor'', Contemp. Math. {\bf 136} (1992), pp.
171--183
\bibitem{esteban} {\it G\'omez Gonz\'alez, E.}, ``Cyclic coverings of a
smooth curve and branch locus of the moduli space of smooth curves'', in
Complex Geometry of Curves, Contemp. Math. {\bf 240} (1999), pp.
183--196

\bibitem{Hejhal} {\it Hejhal, D.A.}, ``Theta functions, kernel functions
and abelian integrals'', Mem. Am. Math. Soc., No. 129 (1972)

\bibitem{KM} {\it Knudsen, F.; Mumford, D.}, ``The projectivity of
the moduli space of stable curves I: preliminaries on {\it det} and
{\it div}'', Math. Scand. {\bf 39} (1976), pp. 19--55

\bibitem{ML} {\it Huybrechts, D.; Lehn, M.}, ``The Geometry of
Moduli Spaces of Sheaves'', Vieweg (1997) Wiesbaden

\bibitem{Laszlo} {\it Laszlo, Y.}, ``Un th\'eor\`eme de Riemann pour les
diviseurs th\^eta sur les espaces de modules de fibr\'es stables sur une
courbe'', Duke Math. J. {\bf 64}  (1991), pp. 333--347

\bibitem{Li} {\it Li, Y.}, ``Spectral curves, theta divisors and Picard
bundles'',
Int. J. Math. {\bf 2}  (1991), pp. 525--550

\bibitem{MP} {\it Mu\~noz Porras, J.M.; Plaza Mart\'{\i}n, F.J.},
``Equations of the moduli space of pointed curves in the infinite
Grassmannian'', J. Differ. Geom. {\bf
51}  (1999), pp. 431--469


\bibitem{Pl} {\it Plaza Mart\'{\i}n, F.J.}, ``Algebraic solutions of the
multicomponent KP hierarchy'', J. Geom. Phys. {\bf 36}
(2000), pp. 1--21

\bibitem{Po} {\it Le Potier, J.},
``Module des fibr\'es semi--stables et fonctions theta'',  in Moduli of
Vector Bundles, Lecture Notes in Pure and App. Math. {\bf
179} (1996), pp. 83--101

\bibitem{Raina-Corr} {\it Raina, A.K.}, ``Fay's Trisecant Identity and
Conformal Field Theory'', Commun.  Math. Phys. {\bf
122} (1989), pp. 625--641

\bibitem{Raina-bc} {\it Raina, A.K.}, ``An Algebraic Geometry Study of
the
$b-c$ system with Arbitrary Twist Fields and Arbitrary Statistics'',
Commun.  Math. Phys. {\bf 140} (1991), pp.
373--397

\bibitem{Raynaud} {\it Raynaud, M.}, ``Sections des Fibr\'es
Vectoriels sur une courbe'', Bull. Soc. Math. France {\bf 110}
(1982), pp. 103--125

\bibitem{Schork} {\it Schork, M.}, ``Generalized $bc$-systems based on
Hermitian vector bundles'', J. Math. Phys. {\bf 41}
(2000),  pp. 2443--2459

\bibitem{Schork2} {\it Schork, M.}, ``On the correlation functions of
the vector bundle generalization of the $b-c$-system'', J. Math. Phys.
{\bf 42} (2001),  pp. 4563--4569

\bibitem{Witten} {\it Witten, E.}, ``Quantum Field Theory,
Grassmannians, and Algebraic Curves'', Commun.  Math. Phys. {\bf 113}
(1988), pp. 529--600



\end{thebibliography}
\end{document}